\documentclass[12pt,reqno,draft]{amsart}

\usepackage{amsmath,amssymb,amscd,latexsym,amsthm}
\usepackage{lineno}
\usepackage{tabularx,euscript}

\theoremstyle{theorem}
\newtheorem{thm}{Theorem}

\theoremstyle{theorem}
\newtheorem{lem}[thm]{Lemma}

\theoremstyle{theorem}
\newtheorem{cor}[thm]{Corollary}

\theoremstyle{remark}

\theoremstyle{remark}
\newtheorem{ex}{Example}

\begin{document}

\title[On a theorem of Morlaye and Joly]{On a theorem of Morlaye and Joly\\ and its generalization}

\author{Ioulia N. Baoulina}

\address{Department of Mathematics, Moscow State Pedagogical University, Krasnoprudnaya str. 14, Moscow 107140, Russia
}

\email{jbaulina@mail.ru}

\date{}

\maketitle

\begin{abstract}
We show that a weaker version of the well-known theorem of Morlaye and Joly on diagonal equations is a simple consequence of a restricted variable version of the Chevalley-Warning theorem. Moreover, we extend the result of Morlaye and Joly to the case of an equation of the form $$b_1D_{m_1}(X_1,a_1)+\dots+b_nD_{m_n}(X_n,a_n)=c,$$ where $D_{m_1}(X_1,a_1)$, \dots, $D_{m_n}(X_n,a_n)$ are Dickson polynomials.
\end{abstract}

\noindent\keywords{\emph{Keywords}:  Equation over a finite field, Chevalley-Warning theorem, diagonal\linebreak equation, Dickson polynomial, value set, combinatorial nullstellensatz}

\noindent\subjclass{\emph{2010 MSC}: 11D79, 11G25, 11T06, 11B75, 12E10, 12E20}

\thispagestyle{empty}

\section{Introduction}

Let $\mathbb F_q$ be a finite field of characteristic~$p$ with $q=p^s$ elements and $\mathbb F_q^*=\mathbb F_q^{}\setminus\{0\}$. In 1934, Artin conjectured that if $H\in\mathbb F_q[X_1,\dots,X_n]$ is a homogeneous polynomial with $\deg(H)<n$ then $H$ has a nontrivial zero. In 1935, Chevalley~\cite{C} proved this and even showed that the hypothesis of homogeneity could be replaced by the weaker assumption of no constant term. Immediately afterwards, Warning~\cite{W} showed that even without this last assumption the characteristic $p$ of $\mathbb F_q$ divides the number of solutions to the equation $H(X_1,\dots,X_n)=0$ in $\mathbb F_q^n$. Both results were extended by their authors to the case of a system of equations of the form $H_i(X_1,\dots,X_n)=0$, $1\le i\le k$, with $\sum_{i=1}^k \deg(H_i)<n$.

In 1999, Alon~\cite{A} established two theorems (so-called \emph{Combinatorial Nullstellensatz}), which provide one of the most powerful algebraic tools, with many applications in combinatorics, number theory and graph theory. Alon used his combinatorial nullstellensatz to give a short proof of the following  weaker version of the Chevalley-Warning theorem: \emph{Let $H_1,\dots,H_k\in\mathbb F_q[X_1,\dots,X_n]$ be polynomials with $\sum_{i=1}^k \deg(H_i)<n$;  if the system of equations $H_i(X_1,\dots,X_n)=0$, $1\le i\le k$, has a solution in $\mathbb F_q^n$, then it has another solution}. Later on, Brink~\cite{Brink} and Schauz~\cite{S} used Alon's ideas to prove a restricted variable version of this theorem.

In 1971, Morlaye~\cite{M} and Joly~\cite{J} established an analogue of the Warning theorem for an equation of the type $b_1X_1^{m_1}+\dots+b_nX_n^{m_n}=c$ (so-called diagonal equation). Namely, they proved that the characteristic $p$ of $\mathbb F_q$ divides the number of solutions to a diagonal equation with exponents $m_1,\dots,m_n$, provided that $1/\gcd(m_1,q-1)+\dots+1/\gcd(m_n,q-1)>1$. One can formulate a weaker version of this result as follows: \emph{Assume that $1/\gcd(m_1,q-1)+\dots+1/\gcd(m_n,q-1)>1$; if a given diagonal equation with exponents $m_1,\dots,m_n$ has a solution in $\mathbb F_q^n$, then it has another solution.}

The goal of this paper is to show that the weaker version of the theorem of Morlaye and Joly stated above and its generalization are consequences of the Chevalley-Warning theorem with restricted variables. In Section~\ref{s3}, we consider the case of a diagonal equation. Section~\ref{s4} contains our main result, namely Theorem~\ref{t2} in which we treat an equation of the type $b_1D_{m_1}(X_1,a_1)+\dots+b_nD_{m_n}(X_n,a_n)=c$, where $D_{m_1}(X_1,a_1)$, \dots, $D_{m_n}(X_n,a_n)$ are Dickson polynomials. Further possible generalizations are discussed in Section~\ref{s5}.

\section{Preliminary lemmas}
\label{s2}

We write $|\EuScript{C}|$ for the cardinality of a finite set $\EuScript{C}$. The following result is a restricted variable version of the Chevalley-Warning theorem (for a proof, see Brink~\cite{Brink} or Schauz~\cite{S}).

\begin{lem}
\label{l1}
Let $H_1,\dots,H_k\in\mathbb F_q[X_1,\dots,X_n]$ be nonzero polynomials and $\EuScript{C}_1,\dots,\EuScript{C}_n$ be nonempty subsets of $\mathbb F_q$ such that
$$
(q-1)\sum_{i=1}^k \deg(H_i) < \sum_{j=1}^n (|\EuScript{C}_j|-1).
$$
Then the set
$$
\{(c_1,\dots,c_n)\in\EuScript{C}_1\times\dots\times\EuScript{C}_n\, |\, H_i(c_1,\dots,c_n)=0\,\,\text{\rm for all $i$}\}
$$
is not a singleton.
\end{lem}

For a polynomial $f\in\mathbb F_q[X]$, let $V_f$ denote the value set of $f$, that is, $V_f=\{f(c)\,|\, c\in\mathbb F_q\}$. The next corollary is a straightforward consequence of Lemma~\ref{l1}.

\begin{cor}
\label{c1}
Let $f_1,\dots,f_n\in\mathbb F_q[X]$ and $H_1,\dots,H_k\in\mathbb F_q[X_1,\dots,X_n]$ be nonzero polynomials. Assume that 
$$
(q-1)\sum_{i=1}^k \deg(H_i) < \sum_{j=1}^n (|V_{f_j}|-1).
$$
If the system of equations
$$
H_i(f_1(X_1),\dots,f_n(X_n))=0,\quad i=1,\dots,k,
$$
has a solution in $\mathbb F_q^n$, then it has another solution. In particular, if $f_1,\dots,f_n$, $H_1,\dots,H_k$ are polynomials without constant terms, then the above system of equations has a nontrivial solution in $\mathbb F_q^n$,  i.e., a solution other than the trivial one $(0,\dots,0)$.
\end{cor}

If $f\in\mathbb F_q[X]$ is a nonconstant polynomial, then each $\gamma\in\mathbb F_q$ has at most $\deg(f)$ preimages under $f$, and so
$$
|V_f|\ge \left\lfloor\frac{q-1}{\deg(f)}\right\rfloor+1,
$$
where $\lfloor z\rfloor$ denotes the greatest integer less than or equal to $z$. Taking $k=1$ and $H_1(X_1,\dots,X_n)=X_1+\dots+X_n$ in Corollary~\ref{c1}, we obtain the following result.

\begin{cor}
\label{c2}
Let $f_1,\dots,f_n\in\mathbb F_q[X]$ be nonconstant polynomials. Assume that
$$
\left\lfloor\frac{q-1}{\deg(f_1)}\right\rfloor+\dots+\left\lfloor\frac{q-1}{\deg(f_n)}\right\rfloor>q-1.
$$
If the equation
$$
f_1(X_1)+\dots+f_n(X_n)=0
$$
has a solution in $\mathbb F_q^n$, then it has another solution. In particular, if $f_1,\dots,f_n$ are polynomials without constant terms, then the above equation has a nontrivial solution in $\mathbb F_q^n$.
\end{cor}

\section{Diagonal equations}
\label{s3}

We consider an equation of the type
\begin{equation}
\label{eq1}
b_1X_1^{m_1}+\dots+b_nX_n^{m_n}=c,
\end{equation}
where $m_1,\dots,m_n$ are positive integers, $b_1,\dots,b_n\in\mathbb F_q^*$, and $c\in\mathbb F_q$. Such an equation is called a \emph{diagonal equation}. Diagonal equations have been studied extensively; see \cite[Chapter~10]{BEW}, \cite[Chapter~6]{LN}, \cite[Chapters~6-7]{MP} and references therein.

Taking $k=1$ and $H_1(X_1,\dots,X_n)=b_1X_1+\dots+b_nX_n-c$ in Corollary~\ref{c1} and using the well-known fact that
\begin{equation}
\label{eq2}
|V_{X^m}|=\frac{q-1}d+1,
\end{equation}
where $d=\gcd(m,q-1)$, we deduce the weaker version of the result of Morlaye~\cite{M} and Joly~\cite{J} mentioned in the introduction.

\begin{thm}
\label{t1}
Assume that
$$
\frac1{d_1}+\dots+\frac1{d_n}>1,
$$
where $d_j=\gcd(m_j,q-1)$, $1\le j\le n$. If the equation~\eqref{eq1} has a solution in $\mathbb F_q^n$, then it has another solution. In particular, if $c=0$, then \eqref{eq1} has a nontrivial solution in~$\mathbb F_q^n$.
\end{thm}

\section{Equations with Dickson polynomials}
\label{s4}

In this section, we consider an equation of the form

\begin{equation}
\label{eq3}
b_1D_{m_1}(X_1,a_1)+\dots+b_nD_{m_n}(X_n,a_n)=c,
\end{equation}
where $m_1,\dots,m_n$ are positive integers, $b_1,\dots,b_n\in\mathbb F_q^*$, $a_1,\dots,a_n,c\in\mathbb F_q$, and $D_{m_1}(X_1,a_1)$, \dots, $D_{m_n}(X_n,a_n)$ are Dickson polynomials defined as follows.

Let $m\ge 1$ be an integer and $a\in\mathbb F_q$. The \emph{Dickson polynomial} 
of degree $m$ with parameter $a$ is defined by
$$
D_m(X,a)=\sum_{j=0}^{\lfloor m/2\rfloor}\frac m{m-j}\binom{m-j}j (-a)^j X^{m-2j}.
$$
Basic properties of Dickson polynomials can be found in the monograph of Lidl, Mullen and Turnwald~\cite{LMT}. Since for $a=0$ we have $D_m(X,a)=X^m$, the equation~\eqref{eq3} can be viewed as a generalization of the diagonal equation~\eqref{eq1}.

Chou, Mullen, and Wassermann~\cite{CMW} used a character sum argument to give bounds for the number of solutions to \eqref{eq3} in $\mathbb F_q^n$. See also \cite{GW} and \cite{OS} for some results on \eqref{eq3} with $a_1=\dots=a_n$, $b_1=\dots=b_n=1$.

In order to generalize Theorem~\ref{t1}, we need the following result concerning the cardinality of the value set of $D_m(X,a)$ with $a\in\mathbb F_q^*$ (for a proof, see \cite[Theorems~10 and $10'$]{CGM} or \cite[Theorems~3.27 and 3.30]{LMT}; see also \cite[Theorem~7]{CMW} for an alternative proof).

\begin{lem}
\label{l2}
Let $a\in\mathbb F_q^*$. Suppose that $2^r\parallel (q^2-1)$. Then
$$
|V_{D_m(X,a)}|=\frac{q-1}{2\gcd(m,q-1)}+\frac{q+1}{2\gcd(m,q+1)}+\delta,
$$
where
$$
\delta=\begin{cases}
1&\text{if $q$ is odd, $2^{r-1}\parallel m$ and $a$ is a nonsquare in $\mathbb F_q$,}\\
1/2&\text{if $q$ is odd, $m$ is even and $2^{r-1}\nmid m$,}\\
0&\text{otherwise.}
\end{cases}
$$
\end{lem}

Combining the relation~\eqref{eq2}, Lemma~\ref{l2} and Corollary~\ref{c1} (with $k=1$ and\linebreak $H_1(X_1,\dots,X_n)=b_1X_1+\dots+b_nX_n-c$), we derive our main result.

\begin{thm}
\label{t2}
Assume that $a_1=\dots=a_t=0$, $a_{t+1},\dots,a_n\in\mathbb F_q^*$, $0\le t\le n$, $2^r\parallel (q^2-1)$ and
$$
\sum_{j=1}^t\frac1{d_j^-}+\frac12\sum_{j=t+1}^n \frac1{d_j^-}+\left(\frac 12+\frac 1{q-1}\right)\sum_{j=t+1}^n\frac1{d_j^+}>1+\frac{n-t-\delta_{t+1}-\dots-\delta_n}{q-1},
$$
where, for $1\le j\le n$, 
\begin{align*}
&d_j^-=\gcd(m_j,q-1),\qquad d_j^+=\gcd(m_j,q+1),\\
&\delta_j=\begin{cases}
1&\text{if $q$ is odd, $2^{r-1}\parallel m_j$ and $a_j$ is a nonsquare in $\mathbb F_q$,}\\
1/2&\text{if $q$ is odd, $m_j$ is even and $2^{r-1}\nmid m_j$,}\\
0&\text{otherwise.}
\end{cases}
\end{align*}
If the equation~\eqref{eq3} has a solution in $\mathbb F_q^n$, then it has another solution. In particular, if $m_{t+1},\dots,m_n$ are odd, $c=0$, and
$$
\sum_{j=1}^t\frac1{d_j^-}+\frac12\sum_{j=t+1}^n \frac1{d_j^-}+\left(\frac 12+\frac 1{q-1}\right)\sum_{j=t+1}^n\frac1{d_j^+}>1+\frac{n-t}{q-1},
$$
then \eqref{eq3} has a nontrivial solution in $\mathbb F_q^n$.
\end{thm}

Note that in some cases the existence of more than one solution can be deduced from Theorem~10 of~\cite{CMW}, while in some other cases the lower bound of Theorem~10 of~\cite{CMW} is trivial, however, our Theorem~\ref{t2} still applies. We illustrate the latter situation with an example.
\begin{ex}
Let $n>4$, and let $e_1,\dots,e_n$ be positive divisors of $q-1$ satisfying $e_1\cdots e_n\le q^{n/2}(q-1)$. For $j=1,2,\dots,n$, put $m_j=(q-1)/e_j$. Let $a_1,\dots,a_n$, $b_1,\dots,b_n\in\mathbb F_q^*$, and put
$$
c=2\sum_{\substack{j=1\\ 2\mid m_j}}^n (-a_j)^{m_j/2}b_j.
$$
Let $N$ denote the number of solutions to the equation~\eqref{eq3} in $\mathbb F_q^n$. Using the notation of Theorem~\ref{t2}, we can rewrite the lower bound of Theorem~10 of~\cite{CMW} as
\begin{equation}
\label{eq5}
N\ge q^{n-1}-q^{(n-2)/2}(q-1)\prod_{j=1}^n (d_j^+ +d_j^-).
\end{equation}
Since $d_j^-=(q-1)/e_j$, we have
$$
d_j^+ +d_j^-\ge 1+\frac{q-1}{e_j}\ge \frac q{e_j}.
$$
Hence
$$
q^{n-1}-q^{(n-2)/2}(q-1)\prod_{j=1}^n (d_j^+ +d_j^-)\le q^{n-1}\left(1-\frac{q^{n/2}(q-1)}{e_1\cdots e_n}\right)\le 0,
$$
and so the inequality~\eqref{eq5} does not give us any information about the number of solutions.

Next, we observe that $d_j^+=\gcd((q-1)/e_j,q+1)\le 2$. Thus
\begin{align*}
\frac 12\sum_{j=1}^n \frac 1{d_j^-}+\left(\frac 12+\frac 1{q-1}\right)\sum_{j=1}^n \frac 1{d_j^+}&\ge \frac{e_1+\dots+e_n}{2(q-1)}+\left(\frac 12+\frac 1{q-1}\right)\cdot \frac n2\\
&\ge \frac n{2(q-1)}+\frac n4+\frac n{2(q-1)}>1+\frac n{q-1}\\
& \ge 1+\frac{n-\delta_1-\dots-\delta_n}{q-1}.
\end{align*}
Note also that \eqref{eq3} has the trivial solution $(0,\dots,0)$. Using Theorem~\ref{t2}, we infer that \eqref{eq3} has a nontrivial solution.
\end{ex}

\section{Concluding remarks}
\label{s5}

It is easy to see that Corollary~\ref{c1} can be extended to the case of restricted variables. For a polynomial $f\!\in\mathbb F_q[X]$ and a nonempty subset $\EuScript{C}$ of $\mathbb F_q$, let
${V_f(\EuScript{C})\!=\!\{f(c)\,|\, c\in\EuScript{C}\}}$. Lemma~\ref{l1} implies

\begin{cor}
\label{c3}
Let $f_1,\dots,f_n\in\mathbb F_q[X]$ and $H_1,\dots,H_k\in\mathbb F_q[X_1,\dots,X_n]$ be nonzero polynomials and $\EuScript{C}_1,\dots,\EuScript{C}_n$ be nonempty subsets of $\mathbb F_q$. Assume that 
$$
(q-1)\sum_{i=1}^k \deg(H_i) < \sum_{j=1}^n (|V_{f_j}(\EuScript{C}_j)|-1).
$$
If the system of equations
$$
H_i(f_1(X_1),\dots,f_n(X_n))=0,\quad i=1,\dots,k,
$$
has a solution with $X_j\in\EuScript{C}_j$, $1\le j\le n$, then it has another solution with $X_j\in\EuScript{C}_j$, $1\le j\le n$.
\end{cor}

Taking $k=1$ and $H_1(X_1,\dots,X_n)=b_1X_1+\dots+b_nX_n-c$ in Corollary~\ref{c3}, we deduce the following extension of Theorem~\ref{t1} to the case where  all the variables are restricted to a subfield of $\mathbb F_q$.

\begin{thm}
\label{t3}
Let $\mathbb F_{p^{\ell}}$ be a subfield of $\mathbb F_q$. Assume that
$$
\frac1{\gcd(m_1,p^{\ell}-1)}+\dots+\frac1{\gcd(m_n,p^{\ell}-1)}>\frac{q-1}{p^{\ell}-1}.
$$
If the equation~\eqref{eq1} has a solution in $\mathbb F_{p^{\ell}}^n$, then it has another solution in $\mathbb F_{p^{\ell}}^n$. In particular, if $c=0$, then \eqref{eq1} has a nontrivial solution in $\mathbb F_{p^{\ell}}^n$.
\end{thm}

It is readily seen that the inequality in Theorem~\ref{t3} can be replaced by
\begin{equation}
\label{eq4}
\frac1{\gcd(m_1,p^{\ell}-1)}+\dots+\frac1{\gcd(m_n,p^{\ell}-1)}>1
\end{equation}
if $b_1,\dots,b_n,c\in\mathbb F_{p^{\ell}}$, but not in general. As an example, we consider an equation of the type $X_1^3+gX_2^3=0$, where $g$ is a generator of the cyclic group $\mathbb F_{25}^*$. This equation has the only solution $(0,0)$ in $\mathbb F_{25}^2$ (and, consequently, in $\mathbb F_5^2$), however, \eqref{eq4} holds for $p^{\ell}=5$, $n=2$, $m_1=m_2=3$.

\section*{Acknowledgments}
The author would like to thank Pete Clark for his interest and encouraging comments. The author thanks the referee for a careful reading of the manuscript and  helpful suggestions.

\end{document}